\documentclass[12pt]{article}
\usepackage{amsxtra,amssymb,amsthm,amsmath,latexsym}
\usepackage{float}
\usepackage{graphicx}
\restylefloat{table}

\theoremstyle{plain}

\def\oH{{\overset{\circ}{H}}}
\def\oH1{{\overset{\circ}{H}\kern-.02in{}^1}}

\def\bee{\begin{equation*}}
\def\eee{\end{equation*}}
\def\be{\begin{equation}}
\def\ee{\end{equation}}

\begin{document}

\title{Numerical results of solving 3D inverse scattering problem with
non-over-determined data
}

\author{Cong Van\\
 Department  of Mathematics, Kansas State University, \\
 Manhattan, KS 66506, USA\\
congvan@math.ksu.edu\\
}

\date{}
\maketitle\thispagestyle{empty}

\begin{abstract}
\noindent We consider the 3D inverse scattering problem with non-over-determined scattering data. The data are the scattering amplitude $A(\beta, \alpha_0, k)$ for all $\beta \in S_\beta^2$, where $S_\beta^2$ is an open subset of the unit sphere $S^2$ in $\mathbb{R}^3$, $\alpha_0 \in S^2$ is fixed, and for all $k \in (a,b), 0 \leq a < b$. The basic uniqueness theorem for this problem belongs to Ramm \cite{R603}. In this paper, a numerical method is given for solving this problem and the numerical results are presented.
\end{abstract}

\section{Introduction}\label{sec1}
Let's first consider the direct scattering problem with a potential:
\begin{eqnarray}
&&(\nabla^2 + k^2 - q(x))u = 0 \quad \text{in} \quad \mathbb{R}^3,\label{eq1}\\
&&u = e^{ik\alpha\cdot x} + v,\label{eq2}\\
&&v = A(\beta, \alpha, k)\frac{e^{ikr}}{r} + o\left(\frac{1}{r} \right), \quad r := |x| \to \infty, \quad \frac{x}{r} = \beta.\label{eq3},
\end{eqnarray}
where $\alpha,\beta \in S^2$ are the directions of the scattered wave and incident wave correspondingly, $S^2$ is the unit sphere, $k^2 > 0$ is energy, $k > 0$ is a constant, $A(\beta,\alpha,k)$ is the scattering amplitude or scattering data, which can be measured, and $q(x) \in Q$, where $Q$ is a set of $C^1$-smooth real-valued compactly supported functions, $q = 0$ for $\max_j|x_| \geq R$, $R > 0$ is a real constant.\\
The direct scattering problem \eqref{eq1} - \eqref{eq3} has a unique solution (see, e.g., \cite{R470}).\\

Consider now the \textit{inverse scattering problem}: find the potential $q(x) \in Q$ from the scattering data $A(\beta,\alpha,k)$.\\
The uniqueness of the inverse scattering problem with fixed-energy data is proved by A.G.Ramm \cite{R470}, i.e., $q(x) \in Q$ is uniquely determined by the scattering data $A(\beta, \alpha, k_0)$ for a fixed $k = k_0 > 0$ and all $\alpha, \beta \in S^2$. A.G.Ramm also gave a method for solving the inverse scattering problem with fixed-energy data and obtained an error estimate for the solution for exact data and also for noisy data, \cite{R584}.\\

\indent In this paper, we give a numerical method for solving the inverse scattering problem with non-over-determined data, i.e., find $q(x) \in Q$ from the scattering data $A(\beta, \alpha_0, k)$ for a fixed $\alpha_0 \in S^2$, all $\beta \in S^2$, and all $k \in (a,b), 0 \leq a < b$. The basic uniqueness theorem for this problem belongs to Ramm \cite{R603}. In Section \ref{sec2}, the idea of the numerical method with its difficulties is presented. In Section \ref{sec3}, the numerical procedure is presented and in Section \ref{sec4}, the numerical results is obtained.

\section{Inversion method}\label{sec2}
Let $D$ be the support of $q(x)$, then the unique solution to \eqref{eq1} - \eqref{eq3} is
\be\label{eq4}
u(x,k) = e^{ik\alpha_0\cdot x} - \int_D g(x,y,k)q(y)u(y,k)dy.
\ee
Let $h(x,k) = q(x)u(x,k)$ then \eqref{eq4} implies
\be\label{eq5}
h(x,k) = q(x)e^{ik\alpha_0\cdot x} - q(x)\int_D g(x,y,k)h(y,k)dy.
\ee
Equations \eqref{eq3} and \eqref{eq4} imply the formula for the scattering amplitude:
\be\label{eq6}
-4\pi A(\beta, k) = \int_D e^{-ik\beta\cdot y} h(y,k)dy.
\ee
Equation \eqref{eq5} implies the formula for finding $q(x)$:
\be\label{eq7}
q(x) = h(x,k)[e^{ik\alpha_0\cdot x} - \int_D g(x,y,k)h(y,k)dy]^{-1}.
\ee
\textit{The idea of our inversion method} is: first discretize \eqref{eq6} to find $h(y,k)$ from the amplitude data $A(\beta, k)$ then use $h(y,k)$ to find $q(x)$ from \eqref{eq7}.\\

Let us partition $D$ into $P$ small cubes with volumn $\Delta_p, 1 \leq p \leq P$. Let $y_p$ is any point inside the small cube $\Delta_p$. Choose $P$ different $k_m \in (a,b), 1 \leq m \leq P$ and choose $P$ different vectors $\beta_j \in S^2m 1 \leq j \leq P$. Then discretize \eqref{eq6} and get
\be\label{eq8}
-4\pi A(\beta_j, k_m) = \sum_{p = 1}^P e^{-ik_m\beta_j\cdot y_p}h_{pm}\Delta_p, \quad 1 \leq j,m \leq P,
\ee
where $h_{pm} = h(y_p, k_m)$. Solve the linear system \eqref{eq8} numerically then use equation \eqref{eq7} to find the values of the unknown potential $q(x_p)$
\be\label{eq9}
q(x_p) = h_{pm}\left[ e^{ik_m\alpha_0\cdot x_p} - \sum_{p' = 1, p' \neq p}^P g(x_p, y_{p'}, k_m)h_{p'm}\Delta_{p'} \right]^{-1}, \quad 1 \leq p \leq P.
\ee
Note that the right hand side of \eqref{eq9} should not depend on $m$ or $j$. This independence is an important requirement in numerical solution of the inverse scattering problem, a compatibility condition for the data. This requirement is automatically satisfied for the limiting integrla equation \eqref{eq7}.\\

The values of $q(y_p)$ essentially determine the unknown potential $q(x)$ if $P$ is large. This potential is unique by the uniqueness theorem in \cite{R603}.\\

Note that one can choose $\beta_j$ and $k_m$ so that the determinant of the system \eqref{eq8} is not equal to zero, so that the system is uniquely solvable, but the difficulty is that the system \eqref{eq8} is very ill-conditioned because it comes from an integral equation of the first kind with an analytic kernel. We use the dynamical system method (DSM) in \cite{R526} to solve the ill-posed system \eqref{eq8}.

\section{Numerical procedure and results}\label{sec3}
In practice, one can measure the scattering data (with noises) experimentally. For our numerical experiments, we need to construct the noisy scattering data $A(\beta_j, k_m)$.
\subsection{Constructing noisy scattering data}\label{sec3.1}
Given a potential $q(x)$, let's first construct the exact scattering data $A^*(\beta_j, k_m)$. We partition $D$ into $P$ small cubes and discretize equation \eqref{eq4} to get
\be\label{eq10}
u(x_p,k_m) = e^{ik_m \alpha_0 \cdot x_p} - \sum_{j = 1}^P g(x_p,y_j,k)q(y_j)u(y_j,k)\Delta_j.
\ee
One solves this linear system to get $u(x_p, k_m)$, then the exact scattering data can be found by the following formula:
\be\label{eq11}
A^*(\beta_j, k_m) = -\frac{1}{4\pi}\sum_{p = 1}^P e^{-ik\beta_j\cdot y_p}q(y_p)u(y_p, k_m)\Delta_p.
\ee
Then one can randomly perturb each $A^*(\beta_j, k_m), 1 \leq j \leq P$ by $\delta^* = $ const $> 0$ to get the noisy scattering data $A(\beta_j, k_m) = A^*(\beta_j, k_m) \pm \delta^*$, here the plus or minus sign is choosen alternatively. So, $||A(\beta_j,k_m) - A^*(\beta_j,k_m)|| = \delta$ for some $\delta > 0$.

\subsection{Numerical procedure}\label{sec3.2}
The following steps are implemented in each experiment
\begin{itemize}
\item[1.] Choose $D$, $\alpha_0$, $P$, $q(x)$, $k_m$ and $\delta^*$.
\item[2.] Use the procedure in Section \ref{sec3.1} to obtain the noisy scattering data $A(\beta_j, k_m)$ with noise $\delta = ||A(\beta_j,k_m) - A^*(\beta_j,k_m)||$.
\item[3.] Try different values of $k_m$ so that the determinant of the system in \eqref{eq8} is not zero. Let $k$ be the found value of $k_m$.
\item[4.] Solve the linear algebraic system \eqref{eq8} to get $h_p = h(y_p, k)$. Here we use the DSM method in \cite{R526} with the noise $\delta$.
\item[5.] Construct the potential $q^*(x_p)$ from the equation \eqref{eq9} with $k_m = k$:
$$q^*(x_p) = h_{p}\left[ e^{ik\alpha_0\cdot x_p} - \sum_{p' = 1, p' \neq p}^P g(x_p, y_{p'}, k)h_{p'}\Delta_{p'} \right]^{-1}, \quad 1 \leq p \leq P.$$
\item[6.] Find the relative error between the reconstructed potential $q^*(x_p)$ and the original potential $q(x)$:
\be\label{eq12}
\text{err} = \frac{||q(x) - q^*(x)||}{||q(x)||} = \sqrt{\frac{\sum_{p = 1}^P|q(x_p) - q^*(x_p)|^2\Delta_p}{\sum_{p = 1}^P|q(x_p)|^2\Delta_p}}.
\ee
\end{itemize}

\section{Numerical results}\label{sec4}
In these experiments, we choose $D$ to be the unit cube around the origin, $P = 1000$, $\alpha_0 = (1,0,0)$, and $50 \leq k_m \leq 100$.
\subsection{For constant potential}\label{sec4.1}
Although the inverse scattering problem with constant potential is not interesting, we use the constant potential to test our inversion method.\\
In this experiment, we take $q(x) = 10$. The following results are obtained:
\begin{table}[H]
\begin{center}\begin{tabular}{ccc}
\hline
\multicolumn{3}{c}{Constant potential $q(x) = 10$}\\
\hline
$\delta^*$ & $\delta  = ||A(\beta_j,k_m) - A^*(\beta_j,k_m)||$ & relative error\\
\hline
0.04 & 0.4348 & 0.0566\\
0.02 & 0.2174 & 0.0037\\
0.01 & 0.1087 & 0.00065\\
\hline
\end{tabular} \caption{Numerical results for constant potential $q(x) = 10$.}\label{tab1}
\end{center}
\end{table}

\begin{figure}[H]
  \includegraphics[width=1\textwidth]{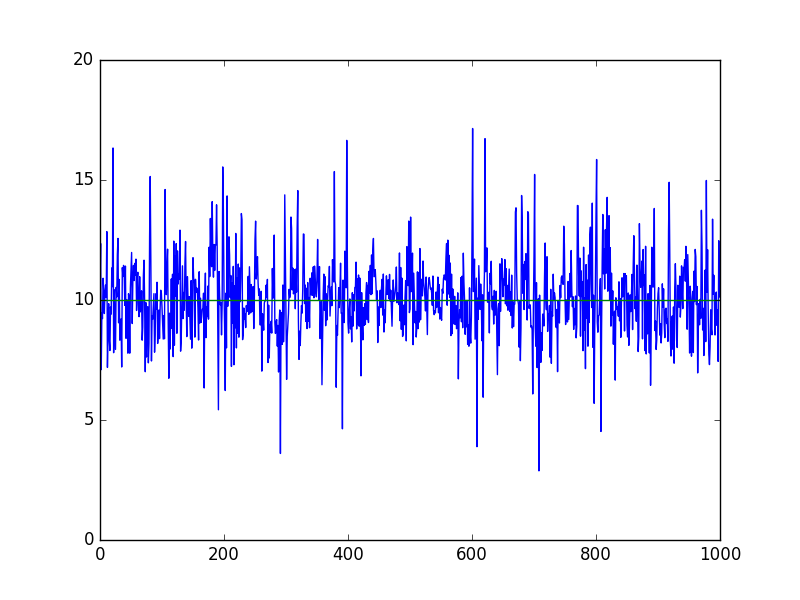}
  \caption{Constructed potential vs original constant potential $q(x) = 10$ when $\delta^* = 0.01$}
\end{figure}

\subsection{For potential $q(x) = \frac{\exp(-|x|)}{|x|}$}\label{sec4.2}
In this experiment, we take $q(x) = \frac{\exp(-|x|)}{|x|}$. The following results are obtained:
\begin{table}[H]
\begin{center}\begin{tabular}{ccc}
\hline
\multicolumn{3}{c}{Potential $q(x) = \frac{\exp(-|x|)}{|x|}$}\\
\hline
$\delta^*$ & $\delta  = ||A(\beta_j,k_m) - A^*(\beta_j,k_m)||$ & relative error\\
\hline
0.04 & 0.0806 & 0.1284\\
0.02 & 0.0403 & 0.0547\\
0.01 & 0.0201 & 0.0367\\
\hline
\end{tabular} \caption{Numerical results for the potential $q(x) = \frac{\exp(-|x|)}{|x|}$.}\label{tab2}
\end{center}
\end{table}

\begin{figure}[H]
  \includegraphics[width=1\textwidth]{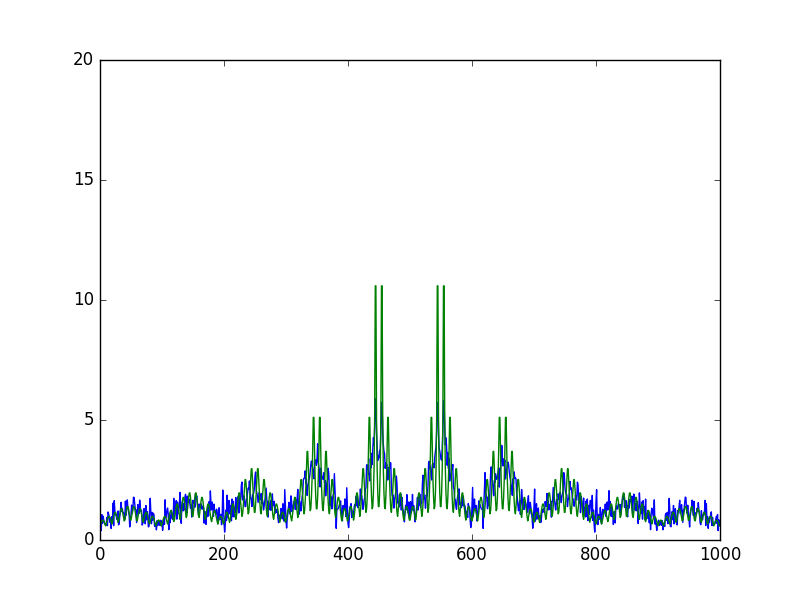}
  \caption{Constructed potential vs original potential $q(x) = \frac{\exp(-|x|)}{|x|}$ when $\delta^* = 0.01$}
\end{figure}

\end{document}